\documentclass[12pt,draft,amsmath,amscd,amsbsy,amssymb]{amsart}
 \relax

\chardef\bslash=`\\ 

\makeatletter \def\verbatim{\interlinepenalty\@M \@verbatim
\leftskip\@totalleftmargin\advance\leftskip2pc \frenchspacing\@vobeyspaces
\@xverbatim} \makeatother \hfuzz1pc

\makeatletter \def\dgt@k{\dg@DX=-3 \dg@DY=2 \dg@SIZE=3} \makeatother

\makeatletter \def\dgt@kk{\dg@DX=3 \dg@DY=-1 \dg@SIZE=3}

\theoremstyle{plain} \newtheorem{thm}{Theorem}[section]
 \newtheorem{lemma}[thm]{Lemma}
\newtheorem{prop}[thm]{Proposition}

\theoremstyle{definition} 
\newtheorem{defin}[thm]{Definition} \newtheorem{que}[thm]{Question}

\newcommand{\R}{{\mathbb R}}

\newcommand{\id}{{\mathrm i}{\mathrm d}}

\newcommand{\U}{{\mathcal U}}

\newcommand{\W}{{\mathcal W}}

\newcommand{\V}{{\mathcal V}}

\numberwithin{equation}{section}

\newcommand{\td}{\tilde}

\newcommand{\CC}{{\mathcal C}}

\newcommand{\mesh}{\mathop{{\mathrm m}{\mathrm e} {\mathrm s}{\mathrm h}}}

\newcommand{\diam}{\mathop{{\mathrm d}{\mathrm i} {\mathrm a}{\mathrm m}}}

\newcommand{\asd}{\mathop{{\mathrm a}{\mathrm s} {\mathrm d}{\mathrm i}
{\mathrm m}}}

\newcommand{\Ind}{\mathop{{\mathrm a}{\mathrm s}{\mathrm I}{\mathrm n}
{\mathrm d}}}

\newcommand{\ind}{\mathop{{\mathrm a}{\mathrm s}{\mathrm i}{\mathrm n}
 {\mathrm d}}}

\newcommand{\N}{{\mathbb N}}

\usepackage[all]{xy}

\usepackage{epic}

\begin{document}


\title[Universal spaces for asymptotic dimension] {Universal spaces for
asymptotic dimension} \author{A.  Dranishnikov} \address{Department of
Mathematics, University of Florida, 358 Little Hall,
PO Box 118105, Gainesville, FL 32611--8105 } \email{dranish@math.ufl.edu}

\author{M.~Zarichnyi}

\address{Department of Mechanics and Mathematics, Lviv National University,
Universytetska 1, 79000 Lviv, Ukraine} \email{topos@franko.lviv.ua,
mzar@litech.lviv.ua}

\subjclass{54F45}

\thanks{The paper was written  during the second author's visit to the
University of
Florida. He thanks the Department of Mathematics  for
hospitality.}
\keywords{Asymptotic dimension, universal space, asymptotic inductive dimension.}

\begin{abstract} We construct a universal  space for the class of proper metric
spaces of bounded geometry and of given asymptotic dimension.  As a consequence
of this result, we establish coincidence of the asymptotic dimension with the
asymptotic inductive dimensions.
\end{abstract}

\maketitle \markboth{A.  Dranishnikov and M.~Zarichnyi} {Universal spaces for
asymptotic dimension}

\section{Introduction}\label{s:intro}

Asymptotic dimension $\asd$ of a metric space was defined by Gromov
for studying asymptotic invariants of discrete groups \cite{G}.  Then a
successful application of the asymptotic dimension was found
by G. Yu \cite{Y}. He proved the Higher Novikov Signature
conjecture for finitely presented groups $\Gamma$ with a finite asymptotic
dimension considered as metric spaces with the word metric. The word metric
$d_S$ on a group $\Gamma$ depends on a generating set $S\subset \Gamma$.
The distance $d_S(x,y)$ in the word metric is the minimal length of
presentations of the word $x^{-1}y$ in the alphabet $S$. It turns out that
the metric spaces $(\Gamma,d_{S})$ and  $(\Gamma,d_{S'})$ are coarsely equivalent
(see Section 2 for the exact definition) for any two finite
generating sets $S$ and $S'$, and hence
$\asd(\Gamma,d_{S})=\asd(\Gamma,d_{S'})$. Thus, in the case of a finitely
generated group $\Gamma$, one can speak about its asymptotic dimension
$\asd\Gamma$ without refereeing to a generating set. In view of Yu's theorem,
a finite-dimensionality results are of a particular interest. In \cite{G} Gromov
proved that all the hyperbolic groups have finite asymptotic dimension.
In \cite{DJ} the finite dimensionality was established for all Coxeter groups.
In \cite{BD1, BD2} finite dimensionality theorems were proved for the fundamental
group of graphs of groups. In particular, it was shown that the asymptotic finite
dimensionality is preserved under the amalgamated product and the HNN-extension.
In \cite{BD2} an upper estimate for $\asd$ of the fundamental group of
a graph of groups is given which turns out to be exact in many cases. Strangely,
it did not give an exact estimate in the case of the free product which is
seemingly the simplest case of the graph  of groups. Only recently the exact
formula has been established in \cite{BDK} with the use of the inductive dimension
$\Ind$ introduced in \cite{D2}. The equality  $\asd=\Ind$ is proven in the present
paper by means of universal spaces.

In the classical dimension theory there are many embedding
theorems of different kind. Here we mension three:

(1) N\"obeling-Pontryagin Theorem: Every compactum $X$ with $\mathrm{dim} X\le n$
can be embedded in $\R^{2n+1}$.

(2) Bowers Theorem \cite{B}: Every compactum $X$ with $\mathrm{dim} X\le n$
can be embedded in the product of $n+1$ dendrites.

(3) Lefschetz-Menger Theorem \cite{Be}: For every $n$ there exists a
universal Menger compactum $\mu^n$, i.e. every  compactum $X$ with $\mathrm{dim} X\le n$
can be embedded in $\mu^n$.

We are interested in the asymptotic analogs of these theorems. The interest is
partly motivated by a topological approach to
Yu's theorem taken in \cite{D1}. First, we note that in the large scale world
there are direct analogs to (1)--(3) for $n=0$. Here we give a description of a space $M^0$
analogous to the Cantor set $\mu^0$. We recall that the Cantor set is the
set of all numbers $0\le x\le1$ which satisfy the following property:

(*) $x$ {\em can be written without use of 1 in the tercimal system
(based on $\{0,1,2\}$)}.

Then $M^0$ can be described as the set of all $x\in\R_+$ satisfying (*).
Every proper metric space $X$ with $\asd X=0$ and of bounded geometry
(see the definition below) can be embedded in the coarse sense in
$M^0$ and hence in $\R$.

The direct analogy between classic and asymtotic dimension
theories ends for $n>0$. It is easy to show that the free group $F_2$ of
two generators has dimension $\asd F_2=1$. It cannot be embedded into $\R^N$,
since $F_2$ has an exponential volume growth. So there is no direct asymptotic
analogy to (1). An asymptotic analogy of (2) was found in \cite{D1}: Every metric
space of bounded geometry can be coarsely embedded in the product of $n+1$
locally finite trees. Since the above product of trees can be embedded in
$(2n+2)$-dimensional nonpositively curved manifold, it gives some hope for
analogy in the case of (1).

The main goal of this paper is an attempt to find a coarse version of the
Lefschetz-Menger theorem, i.e. to construct universal spaces for the asymptotic
dimension. We obtain a rather partial result (Theorem \ref{t:univ2}): For every $n$ there
is a separable metric space $M_n$ with $\asd M_n=n$ which is universal for
proper metric spaces $X$ with bounded geometry and with $\asd X\le n$.
Our universal space is neither proper nor of bounded geometry.

The construction of the space $M_n$ is analogous to those for the
Menger and N\"obeling spaces. Besides the asymptotic strategy, the main
difference is that we are building our fractals not in $\R^{2n+1}$ as in the classical
case but in the product of trees. If one accomplish the local (classical)
construction in the product of (finite) trees or even dendrites he will get
the Menger space according to Bestvina's criterion \cite{Be}. We construct
our spaces $M_n$ out of product of locally finite trees. For
locally finite trees this construction works even better but we have
problems with proving universality. Anyway, we construct for every proper metric space
 of bounded geometry $X$ with $\asd$ an embedding
$X\hookrightarrow\prod_{i=0}^nT_i$ into the product of $n+1$ finite tree
in such a way that $X$ lies in the `$n$-skeleton' $M(\{T_i\})$ of
$\prod_{i=0}^nT_i$. This space   $M(\{T_i\})$ has the property
 $\asd M(\{T_i\})=n$. Moreover, $ \Ind M(\{T_i\})=n$. This allows to prove
the inequality $\Ind X\le\asd X$ for proper metric spaces. For every $X$ the
constructed space     $M(\{T_i\})$ has the Higson property (see Section
\ref{s:Hig}). There were hopes that the mentioned result
on large scale embeddings of asymptotically $n$-dimensional spaces of bounded
geometry into products of $n+1$ nonpositively curved surfaces can be strengthened
so that one can take the hyperbolic planes as these manifolds. It turns out,
however, (see Theorem \ref{t:univ3}) that there is no proper metric space of
finite asymptotic dimension into which every  asymptotically $1$-dimensional
space of bounded geometry can be embedded in the large scale sense.

\section{Preliminaries}\label{s:prelim}

The generic metric will be denoted by $d$. If $A$ is a subset of a metric space
$X$ and $r\in\R$,
then $N_r(A)$ is defined as $\{x\in X\mid d(x,A)<r\}$ if $r>0$ and $A\setminus
\{x\in X\mid d(x,A)<-r\}$ otherwise.  The closed $r$-ball centered at $x\in X$
will be denoted as $B_r(x)$.

A map $f\colon X\to Y$ of metric is called {\em $(\lambda,s)$-Lipschitz} if
$d(f(x),f(y))\le\lambda d(x,y)+s$ for every $x,y\in X$. The
$(\lambda,0)$-Lipschitz maps are also called  $\lambda$-Lipschitz; the 1-Lipschitz maps
are called {\em Lipschitz} or {\em short}. An
{\em asymptotically Lipschitz} map is a map which is $(\lambda,s)$-Lipschitz
for some $\lambda>0$, $s>0$.

For a cover $\U$ of a metric space $X$, we denote by $L(\U)$ its Lebesgue number,
$L(\U)=\inf\{\sup\{d(x,X\setminus U)\mid U\in\U\}\mid x\in X\}$,
and
for a family $\U$ of subsets of a metric space we denote by $\mesh(\U)$ the
lowest upper bound of the diameters of the elements of $\U$.
A family $\mathcal A$ of subsets  of a metric space is called
{\it uniformly bounded\/} if there exists a number $C>0$ such that $\diam (A)<C$
for every $A\in \mathcal A$.

For $r>0$, a metric space  $X$ is called $r$-{\it
discrete\/} if $d(a,b)\ge r$ for every $a,b\in X$, $a\neq b$. A  metric space
$X$ is called {\it discrete\/} if $X$ is $r$-discrete for some $r>0$.

For $r>0$, the $r$-{\it capacity\/} of a subset $Y$ of a metric space is the
maximal
cardinality $K_r(Y)$ of an $r$-discrete subset of $Y$. A metric space $X$ is of
{\it bounded geometry\/} if there exists a number  $r>0$ and a function
$c\colon [0,\infty)\to[0,\infty)$ such that the $r$-capacity of every
$\varepsilon$-ball $B_\varepsilon(x)$ does not exceed $c(\varepsilon)$.

A map $f\colon X\to Y$ between metric spaces is called {\em uniformly
cobounded\/} if for every $r>0$ there exists $C>0$ such that for every $y\in Y$
the diameter of the set $B_r(y)$ does not exceed $C$.

A space $X$ is said to be a {\em geodesic metric space\/} if for
every two points $x,y\in X$ there is an isometric embedding
$j\colon[0,d(x,y)]\to
X$ such that $j(0)=x$ and $j(d(x,y))=y$.

A geodesic segment that connects points $x,y$ in a geodesic metric space will
be denoted by $[x,y]$.

A metric space is called {\em uniformly arcwise connected} if for every $\varepsilon>0$ there
exists $\delta>0$ such that for every two points $x,y\in X$ with $d(x,y)<\varepsilon$ there
exists a path in $X$ connecting $x$ and $y$ and of diameter $\le\delta$.

The (finite) products of metric spaces are always endowed with the
$\sup$-metric.

Let $\CC$ a decomposition of a proper metric space $X$.  Define the {\em
quotient pseudo-metric\/} $\varrho$ on $X$ by the following rule:
$\varrho(x,y)$ is the greatest lower bound of the sums of the form
$\sum_{i=1}^{k-1}d(x_{2i},x_{2i+1})$, where $x=x_1$, $y=x_{2k}$, and for every
$i=1,\dots, k$ there exists an element $C_i\in\CC$ such that
$\{x_{2i-1},x_{2i}\}\subset C_i$.  Obviously, the identity map $\id\colon
(X,d)\to(X,\varrho)$ is short.

A map $f\colon X\to Y$ between metric spaces is called {\em coarse uniform\/} if
for
every $C>0$ there is $K>0$ such that for every $x,x'\in X$ with $d(x,x')<C$ we
have $d(f(x),f(x'))<K$. A map $f\colon X\to Y$  is called {\em metric proper\/}
if the preimage $f^{-1}(B)$ is bounded for every bounded set $B\subset Y$. A map is
{\em coarse} if it is both coarse uniform and metric proper.
Two maps $f,g\colon X\to Y$ between metric spaces are called {\em bornotopic\/}
(see \cite{R}) if there is $C>0$ such that $d(f(x),g(x))<C$ for every $x,y\in
X$. Two metric space $X$ and $Y$ are {\em coarse equivalent\/} if there are
coarse maps
$f\colon X\to Y$, $g\colon Y\to X$ such that the compositions $fg$ and $gf$ are
bornotopic to the corresponding identity maps.

\begin{lemma}\label{l:geod} \cite{D} Let $X$ be a geodesic metric space and $f\colon X\to Y$
be a coarse uniform map. Then $f$ is an asymptotically Lipschitz map.
\end{lemma}

We will need the following simple result whose proof mimicks that of Proposition
2
from \cite{D2}.
\begin{lemma}\label{l:discr} For every proper metric space $X$ and every $C>0$ there is a
$C$-discrete
subset $Y$ of $X$ such that $d(x,Y)\le C$ for every $x\in X$.
\end{lemma}

A map $f\colon X\to Y$ of metric spaces is called a {\em large scale
embedding\/} (see \cite{R}) if there exist increasing functions
$\varphi_1,\varphi_2\colon
[0,\infty)\to[0,\infty)$ with $\lim_{t\to\infty}\varphi_1(t)=
\lim_{t\to\infty}\varphi_2(t)=\infty$ such that $\varphi_1(d(x,y))\le
d(f(x),f(y))\le\varphi_2(d(x,y))$ for every $x,y\in X$. It is easy to see that
every metric space is coarsely equivalent to its image under a coarse embedding.

\subsection{Asymptotic dimension}
The notion of asymptotic dimension is introduced by Gromov \cite{G}.
\begin{defin}
The {\it asymptotic dimension\/} of a metric space $X$ does not exceed $n$
(written
$\asd X\le n$) if for every $D>0$ there exists a uniformly bounded cover $\U$ of
$X$
such that $\U=\U^0\cup\dots\cup\U^n$, where all $\U^i$ are $D$-disjoint.
\end{defin}

A {\em piecewise-euclidean} $n$-dimensional complex of mesh $D$ is a complex which is
a union of   isometric copies of  the standard $n$-dimensional symplex of
mesh $D$ in $\R^{n+1}$.

\begin{lemma}\label{l:g} For a proper metric space $X$ the following two
conditions are equivalent:  \begin{enumerate} \item[1)] $\asd X\le n$; \item[2)]
for any $D>0$ there is a uniformly cobounded short proper map of $X$ into a
piecewise-euclidean complex of mesh $D$ and dimension $n$.  \end{enumerate}
\end{lemma} \begin{proof} See \cite{G}.  \end{proof}

A family of  metric spaces $X_\alpha$ satisfies the inequality $\asd X_\alpha\le n$
{\em uniformly\/} (see \cite{BD1}) if for arbitrary large $D>0$ there exists
$R>0$ and $R$-bounded $D$-disjoint families $\U^0_\alpha,\dots\U^n_\alpha$ such
that $\cup_{i=0}^n\U_\alpha^i$ is a cover of $X_\alpha$.

\begin{thm}\label{t:BD} \cite{BD1} Assume that $X=\cup_\alpha X_\alpha$ and $\asd
X_\alpha$
uniformly. Suppose that for any $R$ there is $Y_R\subset X$ with $\asd Y_R\le n$
and such that the family $\{X_\alpha\setminus Y_R\}$ is $R$-disjoint. Then $\asd
X\le n$.
\end{thm}

\section{Asymptotic embeddings into product of trees}

\subsection{Trees}

A geodesic metric space $T$ is called an {\em $\R$-tree} (a real tree) if (1)
every two points in $T$ are connected by a unique geodesic segment and (2) if
$[x,y]\cap[y,z]=\{y\}$ then $[x,y]\cup[y,z]=[x,z]$ for all $x,y,z\in T$.

Recall that in this definition $[a,b]$ stands for a geodesic segment connecting $a$
and $b$. We will also use self-explaining notations $[a,b)$ and $(a,b)$ for (half)open
geodesic segments.

Every tree (connected connected acyclic graph) is an $\R$-tree. We assume that the graphs are
endowed with a geodesic metric whose restriction to every  edge  is isometric to the standard
unit segment.

A (half)open segment in $T$ is {\em free} if it is an open  subset of $T$.

 The {\em mesh\/} of $T$ is the greatest lower bound of
the diameters of the maximal (with respect to the inclusion) free (half)open segments in $T$.

\begin{defin}\label{d:reg} An  $\R$-tree $T$ is said to be {\em regular\/} if there exists a
sequence of $\R$-trees $T=T_0\supset T_1\supset T_2\supset\dots$ with the following properties:
\begin{itemize} \item[1)] $\cap_{i=0}^\infty
T_i=\emptyset$; \item[2)] for every $i$, the set $T_i\setminus T_{i+1}$ is a
disjoint union of a uniformly bounded family of maximal free in $T_i$ half-open segments;
\item[3)]
$D_0< D_1< D_2<\dots$ and $\lim_{i\to\infty}D_i=\infty$, where $D_i$ denotes the
mesh of the $\R$-tree $T_i$.  \end{itemize} \end{defin}

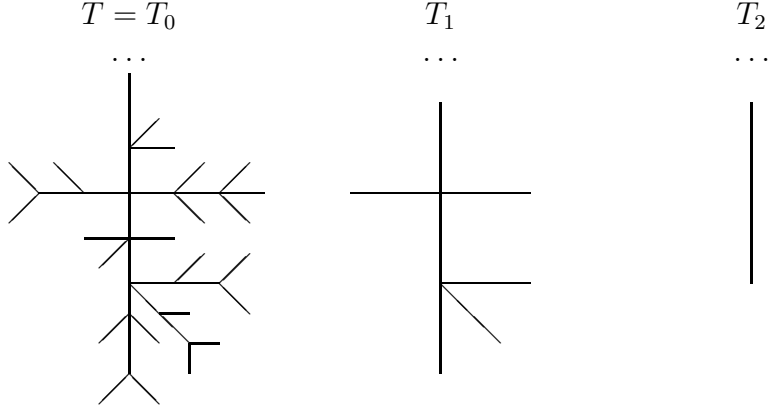
\setlength\unitlength{0.2cm}

\begin{picture}(20,30)
\put(10,23){\makebox(0,5){$\begin{array}{c}T=T_0\\ \dots
\end{array}$}}\drawline(10,23)(10,3)
     \drawline(10,18)(12,20)  \drawline(10,18)(13,18)
  \drawline(4,15)(19,15)  \drawline(4,15)(2,17)  \drawline(4,15)(2,13)
   \drawline(7,15)(5,17)  \drawline(16,15)(18,17)  \drawline(16,15)(18,13)
   \drawline(13,15)(15,17)  \drawline(13,15)(15,13)
   \drawline(7,12)(13,12) \drawline(10,12)(8,10) \drawline(10,9)(16,9)
   \drawline(16,9)(18,11) \drawline(16,9)(18,7)
   \drawline(13,9)(15,11)
   \drawline(10,9)(14,5) \drawline(14,5)(16,5) \drawline(14,5)(14,3)
   \drawline(10,7)(12,5)\drawline(10,7)(8,5)\drawline(12,7)(14,7)
   \drawline(10,3)(12,1) \drawline(10,3)(8,1)
      \end{picture}
\begin{picture}(20,30)
\put(10,23){\makebox(0,5){$\begin{array}{c}T_1\\ \dots
\end{array}$}}\drawline(10,21)(10,3)
       \drawline(4,15)(16,15)  \drawline(10,9)(16,9)
      \drawline(10,9)(14,5)
      \end{picture}
      \begin{picture}(20,22)
\put(10,23){\makebox(0,5){$\begin{array}{c}T_2\\ \dots
\end{array}$}}\drawline(10,21)(10,9)
      \end{picture}

If in the above definition $D_i=2^i$, we say that $T$ is {\em binary regular}.

Suppose that $T^{(0)},\dots,T^{(n)}$ is a sequence of regular $\R$-trees.  For every
$j$,
a sequence of $\R$-trees $T^{(j)}=T_0^{(j)}\supset T_1^{(j)}\supset
T_2^{(j)}\supset\dots$ is given such that properties 1)--3) from Definition
\ref{d:reg} are satisfied with $T_i=T_i^{(j)}$.  For every $i$ and $j$ denote by
$r_i^{(j)}\colon T^{(j)}\to T_i^{(j)}$ the retraction that maps every component
$C$ of the set $T^{(j)}\setminus T_i^{(j)}$ onto its boundary $\partial C$
(which is a singleton). We will refer to these retractions as the canonical
retractions.  For $k\le i$ we have $r_i^{(j)}=r_{ki}^{(j)}r_k^{(j)}$. Let
$$M_j=
\left\{(x_0,\dots,x_n)\in \prod_{i=0}^n
T^{(i)}\mid \text{there exists } i \text{ such that }r_i(x_i)\in\partial
(T^{(i)}_{j-1}\setminus T^{(i)}_j)\right\}$$ and
$M=M(T^{(0)},\dots,T^{(n)})=\cap_{j=1}^\infty M_j$.

Denote by $\mathcal{M}_n$ the class of spaces of the form
$M(T^{(0)},\dots,T^{(n)})$ for some sequence $T^{(0)},\dots,T^{(n)}$ of regular
$\R$-trees.

Suppose  that $T$, $S$ are regular $\R$-trees given with filtrations
$(T_j)_{j=0}^\infty$, $(S_j)_{j=0}^\infty$ satisfying
the conditions from Definition \ref{d:reg}. We denote by $r_j\colon T\to
T_j$, $\varrho_j\colon S\to S_j$ the canonical retractions.

\begin{defin} A map $f\colon T\to S$ is said to be {\em regular\/}
if the following conditions are
satisfied:
\begin{enumerate}
\item[(i)] $f(T_j)\subset S_j$ for every $j$;
\item[(ii)] $f(r_j(\partial (T_{j-1}\setminus T_{j}))\subset
\varrho_j(\partial (S_{j-1}\setminus S_{j}))$  for every $j$.
\end{enumerate}
\end{defin}

The construction $M(T^{(0)},\dots,T^{(n)})$ is functorial in the following
sense. Suppose that $T^i$, $S^i$ are regular $\R$-trees given with filtrations
$(T^i_j)_{j=0}^\infty$, $(S^i_j)_{j=0}^\infty$, $i\in\{0,1,\dots,n\}$ satisfying
the conditions from Definition \ref{d:reg}. We denote by $r_j^i\colon T^i\to
T^i_j$, $\varrho_j^i\colon S^i\to S^i_j$ the canonical retractions. Suppose
that $f_i\colon T^i\to S^i$ are regular maps.
\begin{lemma}\label{l:funct}
Under these conditions $$(f_0\times\dots\times
f_n)(M(T^{(0)},\dots,T^{(n)}))\subset M(S^{(0)},\dots,S^{(n)}).$$
\end{lemma}
\begin{proof} Obvious.
\end{proof}

\begin{prop} For every $M\in \mathcal{M}_n$ we have $\asd M\le n$.  \end{prop}
\begin{proof} We apply Lemma \ref{l:g}.  Note that for every $j$ the map
$$r_j=\prod_{i=0}^nr_j^{(i)}\colon \prod_{i=0}^n T^{(i)}\to \prod_{i=0}^n
T^{(i)}$$ maps $M$ onto an $n$-dimensional piecewise-euclidean polyhedron of
mesh $D_j$.  Since $r_j$ is uniformly cobounded, we are done.  \end{proof}

\begin{thm}\label{t:univ} For every proper metric space $X$ of $\asd X\le n$ there exist
locally finite binary regular trees $T^0,\dots, T^n$ such that $X$ is large scale embeddable
into $M(T^0,\dots,T^n)$.
\end{thm}

The proof is a modification of the corresponding result from \cite{D}.

\begin{prop}\label{p:1} Let $X$ be a proper metric space with base point $x_0$
and \break $\asd X\le n$.  Then there exists a sequence $(\U_k)_{k=1}^\infty$ of open
covers of $X$ such that every $\U_k$ splits into the union
$\U_k=\U_k^0\cup\dots\cup\U_k^n$ of $d_k$-disjoint families and the following
conditions are satisfied:

 \begin{enumerate} \item[1)] $L(U_k)>d_k$ and for every $i$ and every $U\in
\U^i_k$
we have
$N_{-d_k}(U)\neq\emptyset$;
\item[2)] $d_{k+1}=2^{2k+2}m_{k}$, where $m_k$ is the mesh of $\U_k$; \item[3)]
for every $m\in\N$ and for every $i$ there is $k\in\N$ and an element $U\in
\U_k^i$ such that $N_m(x_0)\subset U$.  \item[4)] for every $k,l$, $k<l$ and
every $U\in\U_k^i$, $V\in U_l^i$, if $U\not\subset V$, then $d(U,V)\ge d_k/2$.
\end{enumerate}\end{prop}

\begin{proof} We start with repeating the construction from the proof of
Proposition 1 from \cite{D1}. We proceed by induction. Let $\hat\U_0$ be an open
cover of $X$ which is
$\hat d_0$-discrete with $\hat d_0>2$. We enumerate the partition
$\hat\U_0=\hat\U_0^0\cup\dots\cup\hat\U_0^n$ in such a way that $d(x_0,
X\setminus U)>\hat d_0$ for some $U\in \hat U_0^0$.

 Assume that families $\hat\U_k$ are constructed for all $k\le l$. Define
$d_{l+1}=2^{l+2}m_l$ and consider a uniformly bounded cover $\bar \U_{l+1}$ with
$L(\bar \U_{l+1})>2d_{l+1}$ and such that there is a splitting $\bar
\U_{l+1}=\bar
\U_{l+1}^0\cup\dots\cup\bar \U_{l+1}^n$, where $\bar \U_{l+1}^i$ are
$d_{l+1}$-disjoint and $N_{-2d_{l+1}}(U)\neq\emptyset$ for all $U\in\bar
\U_{l+1}$. The families $\bar \U_{l+1}^0,\dots,\bar \U_{l+1}^n$ are enumerated
in
such a way that $d(x_0,X\setminus U)>2d_{l+1}$ for some $U\in\bar
\U_{l+1}^i$, where $i=l+1\mod n+1$. For every $U\in \bar
\U_{l+1}^i$ let $$\td U=U\setminus \{\overline{N_4(V)}\mid V\in \U^i_k,\ k\le
l,\
V\not\subset U\}$$
and $\U^i_{l+1}=\{\td U\mid U\in\bar
\U_{l+1}^i\}$.

It is proved in \cite{D1} that the sequence $(\hat\U_k)$ satisfies the following
properties:

 \begin{enumerate}
\item[1)] $L(\hat\U_k)>d_k$ and for every $i$ and every $U\in \hat\U^i_k$ we
have
$N_{-d_k}(U)\neq\emptyset$; \item[2)] $d_{k+1}=2^{k}m_{k}$, where $m_k$ is the
mesh of $\hat\U_k$; \item[3)] for every $m\in\N$ and for every $i$ there is
$k\in\N$
and an element $U\in \hat\U_k^i$ such that $N_m(x_0)\subset U$; \item[4)] for
every
$U\in \hat\U_k^i$, $V\in\hat\U^i_l$, $U\not\subset V$ we have $d(U,V)\ge4$.
\end{enumerate}

We are going to modify the sequence $(\hat\U_k)_{k=1}^\infty$.  Passing, if
necessary, to a
subsequence of the sequence $(\hat\U_k)_{k=1}^\infty$ we may assume that the
sequence $(\hat\U_k)_{k=1}^\infty$ itself
satisfies conditions 1), 3),4)  and the following
condition:

$2^*$) $d_{k+1}\ge2^{4k+6}m_{k}$.

Define, by induction, numbers $\td d_k$ and families $\td\U_k^i$.  Let $\td
d_0=d_0$ and $\td\U_0=\U_0$ and suppose $\td d_j$ and $\td\U_j$ have been
defined
for every $j<k$.  For every $U\in\U_k^i$ apply induction by $p$ to define the
sets $\td U(p)$.  Let $\td\U(0)=U$ and suppose sets $\td U(p)$ have been defined
for all $p<q$, for some $q\le k$.  Let $$\td\U(q)=\td U(q-1)\cup\bigcup
\{V\in\td \U^i_{q}\mid d(\td U(q-1),V)<\td d_{q}/2\}.$$ Then, by the definition,
$\td U=\td U(k-1)$ and $\td U_k^i=\{\td U\mid U\in \U_k^i\}$.

We are going to verify conditions
1)--4) (with $U_k$, $d_k$, and $m_k$ replaced by  $\td U_k$, $\td d_k$, and $\td
m_k$ respectively). Conditions 1) and 3) are obvious.

Denote by $\td m_k$ the mesh of $\td \U^i_k$ and let
\begin{equation}\label{f:0}
\td d_k=d_k-2^{k+2}\td m_{k-1}.
\end{equation}
 Note that for every $\td U,\td V\in \td
U_k^i$,
$\td U\neq\td V$ we have
\begin{align*} d(\td U,\td V)=& d(\td U(k-1),\td V(k-1))\ge
d( \td U(k-2),\td V(k-2)) -\td d_{k-1} - 2\td m_{k-1} \\
 \ge&   d(\td U(k-2),\td V(k-2)) - 3\td m_{k-1} \\ \ge& d( \td U(k-3),\td V(k-3)
)
 - 3\td m_{k-1}- 3\td m_{k-2}\\
 & \cdots\\
 \ge &d(U,V)-3\sum_{i=0}^{k-1}\td m_i\ge d_k-3k\td m_{k-1}\ge \td d_k,
\end{align*}
and therefore the family $\td\U_k^i$ is $\td d_k$-disjoint.

Condition 2). Let $\td U\in\td\U^i_k$, then $\td U=\td U(k-1)$ and
\begin{align*} \diam\td U\le &  \diam\td U(k-2)+2\td m_{k-1}+\td d_{k-1}\le
\diam\td
U(k-2)+3\td m_{k-1}\\
&\cdots\cdots\\
\le &\diam \td U(0)+3\sum_{j=0}^{k-1}\td m_j\le m_k+3k\td m_{k-1}\le
m_k+2^{k+2}\td
m_{k-1}
\end{align*}
(here we used an obvious equality $m_j\ge d_j$, for every $j$),
whence
\begin{equation}\label{f:1}
\td m_k\le m_k+2^{k+2}\td m_{k-1}.
\end{equation}
 Then $$\td m_k\le
m_k\sum_{j=0}^{k-1}2^{j(k+2)}\le m_k2^{2k+3}$$
and using condition $2^*$) we obtain
\begin{align}\label{f:2}
\td d_k=& d_k-2^{k+2}\td m_{k-1}\ge d_k-2^{2k+3}m_{k-1} \\ \ge &
2^{4k+6}m_{k-1}-2^{2k+3}m_{k-1} \ge   2^{4k+5}m_{k-1}\nonumber\\ \ge &
2^{2k+2}\td
m_{k-1}.\nonumber
\end{align}

Let us verify condition 4). Suppose the contrary. Then there exist
 $\td U\in\td\U_k^i$, $\td V\in \td \U_l^i$,
$\td U\not\subset \td V$, and $k\le l$ such that $d(\td U,\td V)<d_k/2$. We also
suppose that
$l$ is the minimal number with that property. Since $d(\td U,\td V)<d_k/2$,
there
exists
$x\in \td V$ such that $d(x,\td U)<d_k/2$. Let $j\ge0$ be the minimal integer
with the property $x\in \td V(j)$. It follows from the definition of
$\td \U^i_l$ that $j\le l-1$. Since $x\in  \td V(j)\setminus \td V(j-1)$, it
follows from the definition of $ \td V(j)$ that there exists $\td W\in \td
\U^i_j$
such that $x\in \td W$. Since $j<l$, this contradicts to the choice of $l$.

\end{proof}

For every element $U\in \U_i$ denote by $\psi(U)$ the minimal (with respect to
inclusion) $V\in \U_i$ such that $U$ is a proper subset of $\psi(U)$.

\begin{prop}\label{p:2} Let $X$ be a proper metric space and $(\U_k)$ a sequence of open
covers
of $X$ satisfying properties 1)--4) from Proposition \ref{p:1}.  Let $U\in\U_k$
and $\varrho$ be the quotient metric on $U$ with respect to the decomposition of
$U$ into singletons and the elements of the family $\V= \psi^{-1}(U)$.  For
every $V\in\V\cap\U_j^i$, $W\in\V\cap\U_l^i$, where $j\le l$, we have
$\varrho(V,W)\ge 2^{2j}$.  \end{prop}

\begin{proof} The assertion is obvious for
$j=0$.  Let $V\in\V\cap\U_j^i$, $W\in\V\cap\U_l^i$, where $1\le j\le l$.  Then
$\varrho(V,W)$ is the greatest lower bound of the sums of the form
$\sum_{i=1}^{p-1}d(x_{2i},x_{2i+1})$, where $x\in V$, $x_{2p}\in W$, and for every
$i=1,\dots, p$ there exists an element $C_i\in\psi^{-1}(U)$ such that
$\{x_{2i-1},x_{2i}\}\subset C_i$.  Without loss of generality, we may assume
that in the sums above $C_i\in \psi^{-1}(U)\setminus \V_j$ for all $j<p$.  Then
$(p-2)m_{j-1}+(p-1)d_{j-1}\ge d_j$ and therefore, by condition 2) of Proposition
\ref{p:1}, $$\displaystyle{p-1\ge
\frac{d_j-d_{j-1}}{m_{j-1}+d_{j-1}}\ge\frac{d_j}{4m_{j-1}}\ge 2^{2j}}$$ and
$\sum_{i=1}^{p-1}d(x_{2i},x_{2i+1})\ge (p-1) d_0\ge 2^{2j}$.  Hence
$\varrho(V,W)\ge 2^{2j}$.  \end{proof}

\begin{prop}\label{p:3} Let $X$, $U$, and $\varrho$ be as in Proposition
\ref{p:2}.  There exists $x\in U$ with $\varrho(x,X\setminus U)\ge 2^{2k}$.
\end{prop}

\begin{proof} By condition 1) of Proposition \ref{p:1}, there exists $x\in U$
with
$d(x,\partial U)\ge d_k$. Arguing like in the proof of Proposition \ref{p:2} we
conclude that
$\varrho(x,X\setminus U)\ge 2^{2k}$.
\end{proof}

Fix a sequence $(\U_k)$ of open covers of $X$ satisfying properties 1)--4) from
Proposition \ref{p:1}.  Given $U\in \U^i_k$ define a Lipschitz map $f_U\colon
\bar U\to I_U=[0,2^{2^k}]$ by the following procedure.  Let $\varrho$ denote the
pseudometric generated by the decomposition of $U$ into singletons and the
elements of the family $\psi^{-1}(U)$ and let $\V_j=(\{V\in \U^i_l\mid j\le
l<k\})\cap\psi^{-1}(U)$.

We are going to define by induction with respect to $j$ maps
$\theta_j,\bar\theta_j\colon \partial U\cup\overline{\bigcup\V_j}\to I_U$.  Define
the map $\theta_1\colon \partial U\cup\overline{\bigcup\V_{k-1}}\to I_U$ by the
formula $\theta_1(x)=2^{-k}\varrho(x,X\setminus U)$.  Let
$\bar\theta_1(x)=2^{k-1}[2^{-k+1}\theta_1(x)]$.  The map $\bar\theta_1$ is a
$2^{-k+1}$-Lipschitz map into $I_U$ defined on the closed subset $\partial
U\cup\overline{\cup\V_{k-1}}$ of $\bar U$.

Suppose the maps $\theta_j,\bar\theta_j$ are defined for all $j$, $j_0<j<k$ and
the maps $\bar\theta_j$ are $2^{-j}$-Lipschitz.  By the theorem on extension of
Lipschitz maps (see \cite{D}) there exists a $2^{-j_0+1}$-Lipschitz extension of
$\bar\theta_{j_0-1}$ onto the set $\partial U\cup\overline{\bigcup\V_{j_0}}$.  Let
$\bar\theta_{i_0}(x)=2^{i_0}[2^{-i_0}\theta_{i_0}(x)]$.

As a final result, we obtain a 1-Lipschitz map $$\bar\theta_0\colon \partial
U\cup\overline{\bigcup\V_{0}}=\partial U\cup\overline{\bigcup\psi^{-1}(U)}\to
I_U.$$
Denote by $f_U$ its 1-Lipschitz extension onto $\bar U$.

We will follow \cite{D1} and construct a locally finite tree $T^i$ by the
following
procedure. Assign to every $U\in\U^i_k$, $k=0,1,\dots$,  a line segment
$I_U=[0,2^k]$. The tree
 $T^i$ will be the quotient space of the disjoint union
 $\sqcup\{I_U\mid U\in\cup_{k=0}^\infty\U^i_k\}$ with respect to the following
 equivalence relation $\sim$: suppose $U\in \U^i_k$, $V\in\U^i_l$, where $k<l$
and
$U\in\psi^{-1}(V)$.Then
$0\in I_U$ is identified with $f_V(\bar{\partial U})\in I_V$. We argue like in
the
proof of Theorem 3 from \cite{D1} to show that the quotient space $\sqcup\{I_U\mid
U\in\cup_{k=0}^\infty\U^i_k\}/\sim$ is a  tree. We are going to show that $T^I$
is locally finite. Assume the contrary and let  $y$ be a vertex of infinite
order
in $T^i$. Then there exists an infinite family $\V\subset \U^i$ such that
 $y=q(0_U)$ for every $U\in \V$. Since  the set $\psi^{-1}(U)$ is finite for
 every $U\in\U$, we conclude that there exists an increasing sequence
$(k(j))_{j=1}^\infty$ of positive integers and sets $U_{k(j)}\in\U_{k(j)}^i$
such that $U_{k(j)}\in\psi^{-1}\U_{k(j+1)}$ and $f_{U_{k(j+1)}}(U_{k(j)})=0$.
By conditions 3) and 4), there exists $j>1$ such that $x_0\in U_{k(j)}$ and
$U_{k(1)}\subset U_{k(j)}$. Then, by condition 3), $U_{k(j+1)}\supset
N_{m_{k(j)}+2^{2k(j)+2}d_k}(x_0)$ and, since $\diam(U_{k(j)})\le m_{k(j)}$, by
condition 3), $d(U_{k(j)}, X\setminus U_{k(j+1)})\ge 2^{2k(j)+2}d_k$ and,
therefore, by Proposition \ref{p:3}, $\varrho  (U_{k(j)}, X\setminus
U_{k(j+1)})\ge 2^{k(j)}$. Hence $f_{U_{k(j+1)}}(U_{k(j)})\ge1\neq0$ and we
obtain a contradiction.

Let $q\colon
\sqcup\{I_U\mid U\in\cup_{k=0}^\infty\U^i_k\}\to T^i$ be the quotient map.
Define
the map $f=f^{i}\colon X\to T^i$ as follows. Let $x\in X$. By condition 3) of
Proposition \ref{p:1}, there exists minimal $k$ such that $x\in U$ for some
$U\in\U^i_k$. Put $f^i(x)=qf_U(x)$. It can be easily seen that $f^i$ is
well-defined
and is Lipschitz.

The diagonal map $f=(f^i)_{i=0}^n\colon X\to \prod_{i=0}^n T^i$ is Lipschitz.
Show
that $f$ is a large-scale embedding. Assuming the contrary we can find
sequences
$(x_l)$ and $(y_l)$ in $X$ such that $d(x_l,y_l)\to\infty$ while
$d(f(x_l),f(y_l))<C$
for some $C>0$. Let $k$ be a positive integer such that $2^k>C$. There exists an
integer $l$ such that $d(x,y)>m_{k}$. There exists $i\in\{0,\dots,n\}$ such that
$x_l\in U$ and $d(x,X\setminus U>d_k$ for some $U\in\U^i_k$. Then $f_U(x)=2^k$,
By Proposition \ref{p:3}. Since $d(x,y)>m_k$, we see that $y\notin U$ and
$f_U(y)=0$. Thus, $d(f(x),f(y))\ge2^k>C$, which gives a contradiction.

Note that every $T^i$ is a regular tree. Indeed, for every $j$ denote by $T^i_j$
the
subspace $q(\sqcup\{I_U\mid U\in\cup_{k=j}^\infty\U^i_k\})$. Let us verify
conditions
1)--3) of Definition \ref{d:reg} (with $T$ and $T_j$ replaced by  $T^i$ and
$T^i_j$
respectively). Condition 1) is obvious. The set $T^i_{j+1}\setminus T^i_j$ is a
disjoint union of isometric copies of the half-open segment $(0,2^{j+1}]$ and is
therefore uniformly bounded. To verify condition 3) note that the mesh of the
tree
$T^i_j$ is $2^j$, by the construction.

The trees $T^0,\dots,T^n$ together with filtrations $T_0^i\supset T^i_1\supset
T^i_2\supset\dots$ determine the subspace $M(T^0,\dots,T^n)$ of $\prod_{i=1}^n
T^i$.
We are going to show that $f(X)\subset M(T^0,\dots,T^n)$.  Suppose $x\in X$. For
every $k=0,1,\dots$ there is $i(k)\in\{0,1,\dots,n\}$ such that $x\in U$ for
some $U \in\U^{i(k)}_k$.  Then $r^{i(k)}_{k+1}(f(x))\in \partial
T^{i(k)}_{k+1}\setminus T^{i(k)}_{k}$, i. e. $f(x)\in M_{k+1}$. Therefore,
$f(x)\in M_{k+1}$ for every $k$, i. e. $f(x)\in  M(T^0,\dots,T^n)$.

\subsection{Universal space}

A class of metric spaces $\mathcal{C}$ is said to be universal for a class of
metric spaces $\mathcal{D}$ if for every $D\in \mathcal{D}$ there is a large
scale embedding of $C$ into some $C\in \mathcal{C}$.

Theorem \ref{t:univ} can be reformulated as follows.

\begin{thm}\label{t:univ1} The class $\mathcal{M}_n$ is universal for the class
of proper metric spaces  of asymptotic dimension $\le n$.
\end{thm}

A  metric space $X$ is said to be {\em universal\/} for a class of
metric spaces $\mathcal{D}$ if for every $D\in \mathcal{D}$ there is a large
scale embedding of $D$ into $X$.

Construct an  $\R$-tree $T$ by the following procedure. Define
inductively a sequence $(T_i)$ of $\R$-trees and retractions $r_i\colon T_i\to T_0$.
Let $T_0=\R_+$. Denote by $T_1$ an $\R$-tree obtained by attaching to every integer
point in $\R_+$ a countable set of isometric copies of the unit segment $[0,1]$.
Denote by $r_1$ a retraction of $T_1$ onto $T_0$ that sends every component $C$
which is constant  on every component of the complement $T_1\setminus T_0$.
Suppose that $\R$-trees $T_i$ and retractions $r_i$ are defined for all $i<j$. Denote
by $S$ the subtree $r_{j-1}^{-1}([0,2^{j}])$ with base point $2^j$. The $\R$-tree
$T_j$ is obtained from $T_{j-1}$ by attaching to every point of the form
$k2^j$, $k\in\N$, of $\R_+\subset T_{j-1}$ a countable family of isometric
copies of $S$ by the base point. Denote by $r_j$ a retraction of $T_j$ onto
$T_0$ that sends every component $C$ which is constant  on every component of
the complement $T_j\setminus T_0$.

Let $T=\cup_{i=1}^\infty T_i$. It is easy to see that $T$ is a regular $\R$-tree that
contains  every locally finite binary regular tree  so that the
inclusion preserves the filtration. By Lemma \ref{l:funct}, we conclude that for every
sequence $T^0,T^1,\dots, T^n$ of  locally finite binary regular trees we have
$$M(T^0,T^1,\dots, T^n)\subset M_n=M(S^0,S^1,\dots, S^n),$$ where $S^i$ is
isomorphic to $T$ for every $i=0,1,\dots,n$.

We therefore obtain the following
\begin{thm}\label{t:univ2} There exists a separable metric space $M_n$ with $\asd M_n=n$
universal
for the class of  proper metric spaces $X$ with $\asd X\le n$.
\end{thm}

For the spaces of  asymptotic dimension zero the result can be improved.

\begin{thm}\label{t:0} There exists a proper metric space $M_0$ of bounded
geometry with $\asd
M_0=0$ universal for the class of  proper metric spaces $X$ of bounded
geometry with $\asd X\le 0$.
\end{thm}

\begin{proof} Let $M_0$ denote the set of integers whoce tercimal expansion
consists only from 0s and 1s. Obviously, $\asd M_0=0$ and $M_0$ is of
bounded geometry. For every $k\in\N$ the set $M_0$ can be represented
as the union of $3^k$-discrete family of the intervals of
length $3^k$ in the set of natural numbers.

Let $X$ be a proper metric space of bounded
geometry with $\asd X\le 0$. Let $x_0\in X$ be a base point.  There exists
a sequence $(\U_k)$ of uniformly bounded covers of $X$ with the following
properties:
\begin{enumerate}
\item[(1)] $\U_k$ refines $U_{k+1}$;
\item[(2)] $\U_k$ is $k$-discrete;
\item[(3)] $\cup\{U\in \cup_k\U_k\mid x_0\in U\}=X$.
\end{enumerate}

For $U\in\U_k$, denote by $\psi(U)$ the unique $V\in \U_{k+1}$
that contains $U$.

Since $X$ is of bounded geometry, for every $k\in \N$ there exists
$C_k>0$ such that $|\psi^{-1}(U)|\le C_k$ for every $U\in\U_{k+1}$.

Using the mentioned decomposition of $M_0$ into the union of intervals,
we can easily construct a sequence of covers $(\V_k)$ satisfying the
properties:
\begin{enumerate}
\item[(4)] $\V_k$ refines $V_{k+1}$;
\item[(5)] $\V_k$ is $k$-discrete;
\item[(6)] $|\psi^{-1}(V)|>C_k$ for every $V\in \V_{k+1}$ (as above, for every $W\in
\V_k$ by $\psi(V)$ we denote the unique element of $\V_{k+1}$ that
contains $V$).
\end{enumerate}

For every $U\in \U_1$ and $V\in\V_1$ denote by $g_{UV}$ arbitrary constant map
from $U$ to $V$. Suppose that for every $i\le k$ and every
$U\in \U_i$ and $V\in\V_i$ a map $g_{UV}\colon U\to V$ is defined.  Given
$U\in \U_{k+1}$, $V\in V_{k+1}$, consider arbitrary injective map
$\alpha \colon\psi^{-1}(U)\to\psi^{-1}(V)$. Define $g_{UV}\colon U\to V$
as follows:
$g_{UV}|W=g_{W\alpha(W)}$ for every $W\in \psi^{-1}(U)$.

Now we are ready to define a map $f\colon X\to M_0$. For every $k\in\N$
denote by $U_k$ the unique element of $\U_k$ that contains $x_0$.
By induction, define a sequence of maps $f_k\colon U_k\to M_0$
with the following properties:
\begin{enumerate}
\item[($1'$)] $f_{k+1}|U_k=f|k$, for every $k$;
\item[($2'$)] $f(U_k)$ is contained in an element of $\U_{k+1}$.
\end{enumerate}
Let
$f_1=g_{U_1V}\colon U_1\to V\subset M_0$, for some $V\in\V_1$. Suppose that
$f_i$ are defined for every $i\le k$. Let $V'\in \V_k$ be such that
$f_k(U_k)\subset V'$. Denote by $\alpha \colon \psi^{-1}(\psi(U_k))\to
\psi^{-1}(\psi(V'))$ an embedding such that $\alpha(U_k)=V'$. Define
$f_{k+1}$ by the conditions $f_{k+1}|U_k=f_k$ and $f_{k+1}|W=g_{W\alpha(W)}$
for every $W\in  \psi^{-1}(\psi(U_k))\setminus\{U_k\}$.

 By condition ($1'$), the sequence of maps $(f_k)$ uniquely determines a map
$f\colon X\to M_0$. Using properties (1)--(6) it is easy to see that $f$
is a coarse embedding.
\end{proof}

\subsection{Nonexistence results}

There is no counterpart of Theorem \ref{t:0} in higher dimensions.

\begin{thm}\label{t:1} There is no proper metric space universal for the class
of  proper metric spaces $Y$ with $\asd Y\le 1$.
\end{thm}

\begin{proof} Suppose the contrary and let $X$ be a proper metric space
universal for the class of  proper metric spaces $Y$ with $\asd Y\le 1$. Let
$x_0\in X$ be a base point. Denote by $\alpha(r)$ the 1-capacity of the ball
$N_r(x_0)$, i.e. the number $K_1(N_r(x_0))$. Note that for every $x\in X$ we
have $$K_1(N_r(x))\le K_1(N_{r+d(x,x_0)}(x_0))\le \alpha(r+d(x,x_0)).$$

Let $T$ be a locally finite tree with a base point $t_0$ and the index function
$\varphi(r)$ which we specify later. Suppose that $f\colon T\to X$ is a coarse
embedding. Since $T$ is a geodesic space, the map $f$ is proper and
$(\lambda,s)$-Lipschitz for some $\lambda>0$, $s\ge0$. There exists $a>0$ such
that $d(x,y)\ge a$ implies $d(f(x),f(y))\ge1$. There exist $\lambda',s'>0$ such
that $f(N_r(t_0))\subset N_{\lambda'r+s'}(x_0)$. Therefore
$$K_a(N_r(t_0)\le K_1( N_{\lambda'r+s'}(x_0))=\alpha(\lambda'r+s').$$

Now suppose that $\varphi\colon \N\to\N$ satisfies the following property: for
every $n\in \N$ there exists $m\in \N$ with $\alpha(nr+n)<\varphi(r)$ for all
$r\ge m$.  Define $T$ as follows Let $T_0=\R_+$. For every $m\in\N$ denote by
$T_m$ a tree obtained by attaching to every integer point $n\in\R_+$, $n\ge m$,
$\varphi(n)$ copies of the segment $[0,m]$ by its endpoints. Let
$T=\cup\{T_n\mid n\in\N\cup\{0\}\}$.

Now $$K_a(N_r(t_0))\le \alpha(\lambda'r+s')\le\alpha(nr+n),$$
where $n\in\N$, $n\ge\max\{\lambda',s'\}$. Let $m\in\N$, $m\ge a$. Then for
every $r\ge 2m+1$ we obtain $$K_a(N_r(t_0))\ge K_m(N_r(t_0))\ge\varphi([r]),$$
a contradiction with the choice of $\varphi$.
\end{proof}

Show that there is no universal proper metric space of given asymptotic dimension
in the class of proper metric spaces of bounded geometry.

We fix natural $n\ge1$.

For natural $k$, denote by $X(1,k)$ the set
$$\{(x_0,\dots,x_n)\in [0,k]^{n+1}\mid x_i\notin\mathbb Z\text{ for at
most one }i=0,\dots,n\}.$$
 For natural $m$, let
$X(m,k)=\{mx\mid x\in X(1,k)\}$.

Denote by $\mathbf V$ the set of covers $\V$  of $X(1,k)$ such that
 $\V=\cup_{i=0}^n\V^i$,
where $\V_i$ are $3$-discrete, $i=0,1,\dots,n$. For $\V\in \mathbf V$,
denote by $\mu_{\V}(k)$ the maximal 1-capacity of $V\in \V$ and let
$\mu_k=\min\{\mu_\V(k)\mid \V\in\mathbf V\}$.

Given $\lambda>0$, we say that a subset $A$ of a metric space $X$ is a $<\lambda$-{\em
component}
of $X$ if $A$ is a maximal (with respect to inclusion) subset of $B$ with respect to the
property that every two points of $A$ can be connected by a $\lambda'$-chain, for some
$\lambda'<\lambda$.

\begin{lemma}\label{l:capa} $\mu_k\to \infty$ as $k\to \infty$.
\end{lemma}

\begin{proof} Assume the opposite, i.e. there exists a constant $S$ such that
$\mu_k\le S$, for all $k$. Given
$\V\in\mathbf V$ with $\V_{i=0}^\infty=\cup\V^i$, where $\V^i$
are  $3$-discrete, $i=0,\dots,n$, for every  $V\in\V$, denote  by $\td V$ the family of
$<3$-connected
components of $V$. Let  $\td \V_i=\cup\{\td V\mid V\in \V_i\}$,  $i=0,\dots,n$, and
$\td \V=\cup_{i=0}^n\td \V_i$. From our assumption it easily follows that there exists
 a constant $C>0$
such that for every $W\in \td \V$
we have $\diam(W)\le C$.

Now, for every $W\in\td \V$ define $\hat W$ as the 1-neighborhood of
$W$ in $[0,k]^{n+1}$. The families
$\{\hat W\mid W\in\td \V_i\}$, $i=0,1,\dots,n$, are 1-discrete, of mesh $\le C$,
and together form a cover of  $[0,k]^{n+1}$. Then the  families
$\{(1/k)\hat W\mid W\in\td \V_i\}$, $i=0,1,\dots,n$, are discrete, of mesh $\le C/k$,
and together form a cover of  $[0,1]^{n+1}$. Because of arbitrarity of $k$, we
obtain
a contradiction with Ostrand's characterization of the covering dimension
\cite{O}.
 \end{proof}

Passing to the images of spaces under the homothety with coefficient $m$ and centered at
the origin, we derive from Lemma \ref{l:capa} the following statement.

\begin{prop}\label{p:capa} For every $m,s\in \N$ there exists $r(m,s)\in \N$ such that for
any
$r\ge r(m,s)$ the following holds: given a cover $\V$ of $X(m,r)$ that splits into union
of $n+1$ $3m$-disjoint families, there is $V\in \V$ with $K_m(V)\ge s$.
\end{prop}

\begin{thm}\label{t:univ3} There is no proper metric
space $Y$ of bounded geometry,  \break $\asd Y=n$, with
the following property: for every proper metric space $X$  of bounded geometry,
$\asd X=n$,
there exists a large scale embedding of $X$ into $Y$.
\end{thm}

\begin{proof} Suppose the opposite and let $Y$ be such a space. There exists a
sequence of covers $(\U_k)$ of $Y$  with the following properties:
\begin{enumerate}
\item[1)] $\U_k$ is uniformly bounded;
\item[2)] $\U_k=\cup_{i=0}^n\U_k^i$, where $\U_k^i$ are $d_k$-discrete with $d_k\ge k^2$.
\end{enumerate}

Since $Y$ is of bounded geometry, $c_k=\max\{K_1(U)\mid U\in\U_k\}<\infty$.

For every $m\in\N$, choose $R_m$ so that  for any
$r\ge r(m,R_m)$ the following holds: given a cover $\V$ of $X(m,r)$ that splits into union
of two $3m$-disjoint families, there is $V\in \V$ with
$K_m(V)\ge \sum_{i=1}^m c_m+1$.

Define a space $X$  as follows: attach to $\R_+$ at every point
$m\in\N\subset\R_+$
a copy of space
$X(m,R_m)$, with the geodesic metric on the adjunction space.
Note that $\asd X=1$.

Suppose that $f\colon X\to Y$ is a large scale embedding, i. e. there exist
monotone, increasing to infinity functions $\varphi_1,\varphi_2\colon\R_+\to\R_+$
such that $$\varphi_1(d(x,y))\le d(f(x),f(y))\le \varphi_2(d(x,y))$$
for every $x,y\in X$.

Since $X$ is a geodesic space, by Lemma \ref{l:geod} one may choose $\varphi_2$ linear, $\varphi_2(t)=at+b$, $a>0$.
There exists $m_0\in \N$ such that $d(f(x),f(y))\ge1$ as $d(x,y)\ge m_0$.

For any $m\ge m_0$, find minimal $k(m)\in \N$ such that the preimage $f^{-1}(\mathcal A)$
of any
$d_{k(m)}$-discrete family is $3m$-discrete. We have for $x\in A$, $y\in B$,
$A,B\in\mathcal A$,
$A\neq B$:
$$d_{k(m)}\le d(f(x),f(y))\le\varphi_2(d(x,y))=ad(x,y)+b$$
whence $d(x,y)\ge (d_{k(m)}-b)/a$ and if $(d_{k(m)}-b)/a\ge 3m$, we are done. Therefore,
the minimal $k(m)$ is found from the inequality $d_{k(m)}\ge 3ma+b$.

We see that there exists $m_1\in\N$ such that
$d_m\ge (3m-b)/a$ for every $m\ge m_1$.

Now consider $m>\max\{m_0,m_1\}$. The cover $\V_m$
splits into the union of  $d_m$-disjoint
families $\V_m^i$, $i=0,1\dots, n$.

The cover $f^{-1}(\V_m)$ of
$X(m,R_m)$ splits into the union of  $3m$-disjoint families
$f^{-1}(\V_m^i)$, $i=0,\dots,n$.
 By Proposition \ref{p:capa}, there exists an element $V\in \V_m$ such that
$K_m(f^{-1}(V))\ge \sum_{i=1}^mc_i+1>c_m$. As $m>m_0$, we see that $K_1(V)>c_m$,
a contradiction.

\end{proof}

\section{Higson property}\label{s:Hig}
A  metric space $X$ with $\asd X\le n$ is said to satisfy the {\it Higson
property\/} if there exists  $C>0$ such that
for every $D>0$ there exists a cover $\U$ of $X$ with $\mesh(\U)<CD$ and
such that $\U=\U^0\cup\dots\cup\U^n$, where  $\U^0,\dots,\U^n$ are $D$-disjoint.
Equivalently, $X$ has the Higson property if there exists a sequence
$(\U_k)$
of uniformly bounded covers of $X$ and $C>0$ such that
$\mesh(\U_k)\to\infty$ and every ball in $X$ of radius
$\le C\mesh(\U_k)$ intersects at most $n+1$ element of $
(\U_k)$. The latter is a large scale analog of the
{\em Nagata dimension}  N-$\dim$
defined by Assouad  \cite{a}
  as follows: for a
metric space $X$ we have N-$\dim X\leq n$ if there
is a constant $C$ such that, for each $r>0$, there is an open cover $\U(r)$ of
$X$ by sets of diameter $\le Cr$ such that
each open ball of radius $r$ meets at most $n+1$ members of $\U(r)$.

\begin{prop}\label{p:hig} Every proper metric space
 of finite
asymptotic
dimension is coarsely equivalent to a proper metric space that satisfies
the Higson
property.
\end{prop}
\begin{proof} By Theorem \ref{t:univ}, there is a large scale embedding of $X$
into a space of the form  $M(T^0,\dots,T^n)$, for  a sequence of regular locally
finite trees $T^0,\dots,T^n$ together with filtrations $T_0^i\supset
T^i_1\supset
T^i_2\supset\dots$. Since the Higson property is preserved by subspaces, it is
sufficient to show that the space  $M(T^0,\dots,T^n)$ satisfies it.

Note that the set $K_j=r_j(M(T^0,\dots,T^n)$ is a cubic piecewise euclidean
$n$-dimensional complex with $\mesh(K_j)=D_j=2^j$. Passing to the second
barycentric
subdivision of $K_j$ we are able to produce $D_k/3$-discrete families $\mathcal
A_l$ with  $\mesh(\mathcal A_l)\le 2D_k/3\le D_k$, $l\in\{0,\dots,n\}$. Let
$\U_k^l=\{r_k^{-1}(A)\mid A\in \mathcal A_l\}$, then
$$\mesh(\U_k^l)\le2(D_k+D_{k-1}+\dots+D_{0})\le4D_k.$$

Since the map $r_k$ is short, the family $\U_k^l$ is also $D_k/3$-discrete.
\end{proof}

A metric space $X$ satisfies the {\it $n$-dimensional Nagata property} (briefly
denoted by $(\mathrm P_n)$) if for every $r>0$, every    $x\in X$, and every
$y_1,\dots,y_{n+2}$ such that $d(y_i,N_r(x))<2r$ for every $i=1,\dots,n+2$, there exist
$i,j\in\{1,\dots,n+2\}$, $i\neq j$, such that $d(y_i,y_j)<2r$ (see \cite{N}).

\begin{thm} For a  proper metric space $X$ the following  are equivalent:
\begin{enumerate}
\item[(1)] $\asd X\le n$;
\item[(2)] $X$ is large scale equivalent to a metric space with property
$(\mathrm P_n)$.
\end{enumerate}
\end{thm}
\begin{proof} (1)$\Longrightarrow$(2) By Proposition \ref{p:hig}, $X$ is coarsely
equivalent to a proper metric space with the Higson property.  Moreover, as we can see
from the proof of  Proposition \ref{p:hig}, without loss of generality,
we may assume that there exists a sequence $(\U_k)$ of uniformly bounded open
covers of $X$ and $C>0$ such that $\mesh(\U_k)=2^k$ and $\U_k$ splits into the union
of $n+1$ families that are $C2^k$-discrete. Equivalently, every ball of radius
$C2^{k-1}$ intersects at most $n+1$ element of the cover $\U_k$.  By
Lemma \ref{l:discr}, there exists a 1-discrete subset $X'$ of $X$ which is coarsely
equivalent to $X$. For every $r\in(0,\infty)$, put $$\U'(r)=\begin{cases}\{\{x\}\mid
x\in X'\}& \text{ whenever }0< r <1,\\
\{U\cap X'\mid U\in\U_{[r]}\}&   \text{ whenever }r\ge1,
\end{cases}$$
and $C'=C/4$. Then every ball in $X'$ of radius $\le C'r$ intersects at most $n+1$
element of the cover $\U'(r)$, $r\in(0,\infty)$. The latter means that
N-$\dim X'\le n$ and by

It is proved in \cite{a}[Proposition 2.2] that then there exists
 a metric  $\delta$ on
$X'$ such that $\delta$ satisfies $(\mathrm P)_n$ and $\delta^p$ is Lipschitz 
equivalent to the original metric, $d$, on $X'$. The latter means that there
exist $A,B>0$ such that 
$$A\delta(x,y)^p\le d(x,y)\le B\delta(x,y)^p, \ x,y\in X'.$$

(2)$\Longrightarrow$(1) The proof is completely parallel to that
 in the  case of dimension $\dim$ (see \cite{N}). Let $X$ be a proper
 metric space whose metric, $d$, satisfies property $(\mathrm P_n)$. Given $D>0$,
 find a maximal (with respect to inclusion) set $Y\subset X$ which is $2D$-discrete.
Then $\U=\{N_{2D}(y)\mid y\in Y\}$ is a uniformly bounded cover  of $X$.

Show that for every $x\in X$ the set $N_D(x)$ intersects at most $n+1$ element
of  $\U$. Indeed, suppose the contrary; then there $y_1,\dots,y_{n+2}\in Y$
such that  $N_D(x)\cap N_{2D}(y_i)\neq\emptyset$ for every $i=1,\dots,n+2$.
Property $(\mathrm P_n)$ implies that there exist $i,j\in\{1,\dots,n+2\}$, $i\neq j$,
such that $d(y_i,y_j)<2D$. This is a contradiction with the choice of $Y$.
Because of arbitrarity of $D$, $\asd X\le n$.

\end{proof}

\section{Coincidence of asymptotic dimensions}

\subsection{Higson compactification and Higson corona}
Let $\varphi\colon X\to\R$ be a function defined on a metric space $X$.
For every $x\in X$ and every $r>0$ let $V_r(x)=\sup\{|\varphi(y)-\varphi(x)|\mid
y\in N_r(x)\}$. A function $\varphi$ is called {\it slowly oscillating\/}
whenever for every $r>0$ we have
$V_r(x)\to0$ as $x\to\infty$ (the latter means that for every $\varepsilon>0$
there exists a compact subspace $K\subset X$ such that $|V_r(x)|<\varepsilon$
for all $x\in X\setminus K$. Let $\bar X$ be the compactification of $X$ that
corresponds
to the family of all continuous bounded slowly oscillation functions. The {\it
Higson
corona\/} of $X$ is the remainder $\nu X=\bar X\setminus X$ of this
compactification.

It is known that the Higson corona is a functor from the category of proper
metric space and coarse maps into the category of compact Hausdorff spaces. In
particular, if $X\subset Y$, then $\nu X\subset \nu Y$.

For any subset $A$ of $X$ we denote by $A'$ its trace on $\nu X$, i. e. the
intersection of the closure of $A$ in $\bar X$ with $\nu X$. Obviously, the set
$A'$ coincides with the Higson corona $\nu A$.

\subsection{Asymptotic inductive dimensions}
The notion of asymptotic inductive dimension $\Ind$ is introduced in \cite{D2}.

Recall (see, e.g. \cite{E}) that a closed subset $C$ of a topological space $X$ is a {\it
separator\/}
between disjoint subsets $A, B\subset X$ if $X\setminus C=U\cup V$, where $U,V$
are open subsets in $X$, $U\cap V=\emptyset$, $A\subset U$, $V\subset B$.

 Let $X$ be a proper metric space. A subset
$W\subset X$ is called an {\it asymptotic neighborhood\/} of a subset $A\subset
X$ if $\lim_{r\to\infty}d(X\setminus N_r(x_0),X\setminus W)=\infty$. Two sets
$A,B$
in a metric space are {\it asymptotically disjoint\/} if $\lim
_{r\to\infty}d(A\setminus N_r(x_0),B\setminus N_r(x_0))=\infty$. In other words,
two
sets are asymptotically disjoint if the traces $A'$, $B'$ on $\nu X$  are
disjoint.

A subset $C$ of a metric space $X$ is an  {\it asymptotic separator\/}  between
asymptotically disjoint subsets $A_1,A_2\subset X$ if the trace $C'$ is a
separator
in $\nu X$ between $A_1'$ and $A_2'$.

 By the
definition, $\Ind X=-1$ if and only if $X$
is bounded. Suppose we have defined the class of all proper metric spaces $Y$
with $\Ind Y\le n-1$. Then $\Ind X\le n$
 if and only if for every asymptotically
disjoint subsets $A_1,A_2\subset X$ there exists an  asymptotic separator
 $C$
between
$A_1$ and $A_2$ with $\Ind C\le n-1$.
The dimension function $\Ind$ is called the {\em asymptotic
inductive dimension\/}.

We can similarly define the small inductive asymptotic dimension $\ind$. By the
definition, $\ind X=-1$ if and only if $X$
is bounded. Suppose we have defined the class of all proper metric spaces $Y$
with $\ind Y\le n-1$. Then $\ind X\le n$
 if and only if for every  $A\subset X$ and every $x\in\nu X\setminus A'$ there
exists a subset $C$ of $X$ such that $C'$ is a separator between $A'$ and
$\{x\}$ in $\nu X$ and $\ind C\le n-1$.

\begin{prop} For every metric space $X$ we have $\ind X\le \Ind X$.
\end{prop}

\begin{proof} Induction by $\Ind X$. By the definition, the properties $\ind
X=-1$ and $\Ind X=-1$ are equivalent. Suppose $\Ind X\le n$ and $A\subset X$,
$x\in \nu X\setminus A'$. There exists a continuous function $\varphi\colon \bar
X\to [-1,1]$ such that $\varphi|\bar A=-1$, $\varphi(x)=1$. The sets $A$ and
$B=(\varphi|X)^{-1}([0,1])=\varphi^{-1}([0,1])\cap X$ are asymptotically
disjoint subsets of $X$ ad therefor there exists an asymptotic separator $C$
between $A$ and $B$ with $\Ind C\le n-1$. Then obviously, $C'$ is a separator
between $A'$ and $x$ in $\nu X$ and, by inductive hypothesis, $\ind C\le n-1$.
Therefore $\ind X\le n$.
\end{proof}

The following is a counterpart of Theorem 1 in \cite{D2}.
\begin{prop} Let $X$ be a proper metric space. Then $\ind X\ge {\mathrm{ind}}\nu
X$.
\end{prop}
\begin{proof} Induction by $\ind X$. Obviously, $\ind X=-1$ if and only if
${\mathrm{ind}}\nu X=-1$. Suppose that $\ind X\le n$ and show that
${\mathrm{ind}}\nu X\le n$. Let $A$ be a closed subset of $\nu X$ and $x\in\nu
X\setminus A$. There exists a continuous function $\varphi\colon \bar X\to
[-1,1]$ such that $\varphi|A=-1$, $\varphi(x)=1$. Let $B=\varphi^{-1}([-1,0])$.
Then $x\notin B'$ and, since $\ind X\le n$, there exists a subset $C\subset X$
with $\ind C\le n-1$ and such that $C'$ is a separator between $B'$ and $x$ in
$\nu X$. By the inductive hypothesis, $${\mathrm{ind}}C'={\mathrm{ind}}\nu C\le
\ind C\le n-1.$$ Since $C$ is a separator between $B'$ and $x$ in $\nu X$, it is
also a separator between $A$ and $x$.

Therefore, for every closed subsert $A$ of $\nu X$ and $x\in \nu X\setminus A$
we are able to find a separator $K$ with ${\mathrm{ind}}K\le n-1$. This means
that ${\mathrm{ind}}\nu X\le n$.
\end{proof}

It is proved in \cite{D2} (see Proposition 1 therein) that $\Ind X=\Ind Y$ for
coarsely equivalent spaces $X$ and $Y$.

\begin{prop} Let $X$ and $Y$ be coarsely equivalent metric spaces. Then $\ind
X=\ind Y$.
\end{prop}
\begin{proof} Analogous to that of Proposition 1 in \cite{D2}.
\end{proof}

\begin{lemma}\label{l:ass} Suppose that $A\cup B$ are asymptotically disjoint
subsets of
$X\subset Y$ and $C$ is an asymptotic  separator
in $Y$ between $A$ and $B$ with $\asd C\le m$. Then there exists an  asymptotic
 separator $\td C$
 between $A$ and $B$ in $X$ with $\asd \td C\le m$
\end{lemma}

\begin{proof} Define subsets $D_k$ of $X$ by induction. Let $D_1=N_1(C)\cap X$
and $D_{k+1}=(N_{k+1}(C)\cap X)\setminus N_k(D_k)$. By Lemma \ref{l:discr}, there
exists a $(k+1)$-net $\td D_{k+1}$ in $D_{k+1}$ which is $k$-discrete. Let $\td
C=\cup_{k=1}^\infty\td D_k$.

We first show that $\asd \td C\le m$. Remark that for every $k>0$ the complement
to the subset $\cup\{\td D_j\mid j\le k\}$ is  the union of a $k$-disjoint
family $\{\td D_j\mid j>k\}$. To apply Theorem \ref{t:BD}, we have to verify that
$\asd \td D_k\le m$ uniformly. If $k>0$, then for some $R>0$ there exists an
$R$-bounded cover $\U$ of the set  $\cup\{\td D_j\mid j\le k\}$ that splits into
the union $\U^0\cup\dots\cup\U^m$ of $k$-discrete families. Let
$\V^0=\U^0\cup\{\{x\}\mid x\in \cup\{\td D_j\mid j>k\}\}$, then the family
$\V^0\cup\U^1\cup\dots\cup\U^m$ is a $\min\{R,k\}$-uniformly bounded family that
splits into the union of $m+1$ $k$-disjoint families.

Show that for every integer $k>0$ there exists $R=R(k)>0$ such that $N_k(C)\cap
X\subset N_{R(k)}(\td C)$. Indeed, if $x\in N_k(C)\cap X$, then either $x\in
D_k\subset N_{2k}(\td D_k)$ or $$x\in N_k(C)\setminus D_k\subset
N_{k-1}(D_{k-1})\subset N_{k-1}(N_{k-1}(\td D_{k-1}))\subset N_{2k-2}(\td
D_{k-1})\subset N_{2k-2}(\td C)$$
and we can choose $R(k)=2k$.

Now suppose that $C$ is an asymptotic  separator between $A$ and $B$ in $Y$.
Show that $\td C$ is an asymptotic  separator between $A$ and $B$ in $X$. It is
sufficient to show that $\td C'\supset C'\cap\nu X$. Suppose that $a\in
C'\cap\nu X$ and $U$ is a closed neighborhood of $a$ in $\bar Y$. Then there
exist sequences $(x_i)$ in $X\cap U$ and $(c_i)$ in $C\cap U$  and $k\in\N$ such
that $d(x_i,c_i)\le k$ for every $i$.

Case 1). For infinite number of $i$ (passing to a subsequense we then assume
that for all $i$) $x_i\in D_k$. Then there exist $\td x_i\in\td D_k$ with
$d(x_i,\td x_i)\le k$.

Case 2). For infinite number of $i$ (passing to a subsequense we then assume
that for all $i$) $x_i\in N_{k-1}(D_{k-1})$. Then there exist $y_i\in D_{k-1}$
such that $d(x_i,y_i)\le k-1$ and there exist $\td x_i\in\td D_{k-1}$ such that
$d(y_i,\td x_i)\le k-1$.

In both cases, $\td x_i\in\td C$ for every $i$ and $d(x_i,\td x_i)\le
\max\{k,2k-1\}=2k$ and therefore $\td x_i\in U\cap X$ for all but finitely many
$i$. This means that the sets $\td C$ and $U$ are not asymptotically disjoint
and therefore $C'\cap U\neq\emptyset$. Because of arbitrarity of $U$ we conclude
that $a\in\td C'$.
\end{proof}

\begin{thm}\label{t:ind} For every space $X\in \mathcal{M}_n$ we have $\Ind X\le
n$.
\end{thm}
We need some auxiliary results. Let $T^0,\dots,T^n$ be a sequence of regular locally finite
trees.

\begin{prop}\label{p:assep} For every disjoint
asymptotically disjoint subsets $A,B$ in $M(T^0,\dots,T^n)$ there exists a
separator $C$ between them such that $\asd C\le n-1$ and the sets $A$, $B$, and
$C$ are asymptotically disjoint.
\end{prop}
\begin{proof}
Let \begin{align*}
U=&\{x\in M(T^0,\dots,T^n)\mid d(x,A)\le (1/2)d(x,B)\},\\ V=&\{x\in
M(T^0,\dots,T^n)\mid d(x,B)\le (1/2)d(x,A)\},\end{align*} then $U,V$ are closed,
disjoint, asymptotically disjoint asymptotic neighborhoods of the sets $A,B$
respectively.  For every $j$ let $R_j$ be a number such that
$$d(r_j(U\setminus N_{R_j}(x_0)),r_j(V\setminus N_{R_j}(x_0)))\ge2^{j+2}$$ (such
a number exists because $U,V$ are asymptotically disjoint and the
retraction $r_j$ is $2^{j+1}$-close to the identity map. We will assume that
$R_{j+1}\ge 2R_j$. Put $U_j=r_j\left(U\cap(\overline{N_{R_{j+1}}(x_0)}\setminus
N_{R_{j}}(x_0))\right)$.

Denote by $\hat U_j$ the union of all $n$-dimensional cubic simplices in the
subcomplex
$r_j(M(T^0,\dots,T^n))$ that intersect $U_j$. Let
$S_j=r^{-1}(\partial(N_{1/2}(\hat U_j))$ and $S=\cup_{i=0}^\infty S_j$.

We are going to show that $S$ contains a separator between $U$ and $V$ and
$\asd(S)\le n-1$.

Put $S'=\cup\{S_{2k}\mid k\in\N\}$, $S''= \cup\{S_{2k-1}\mid k\in\N\}$.
By finite sum theorem for asymptotic dimension (see, e.g. \cite{BD1}), it is
sufficient to
prove that $\asd S'\le n-1$, $\asd S''\le n-1$ and, because of complete analogy,
we only prove that $\asd S'\le n-1$.

Prove that the family $\{S_{2k}\}$ has asymptotic dimension $\le n-1$ uniformly.
Let $D>1$. There exists natural $j$ such that every $(n-1)$-dimensional cubic
complex of mesh $2^j$ can be covered by a $2^j$-bounded family of its subsets
that splits into the union of $n$ families which are $2D$-disjoint.

{\bf Claim.\/} $r_{jj'}^{-1}(\partial(N_{1/2}(\hat U_{j'}))\subset
N_{3/4}(r_{jj'}^{-1}(\hat
U_{j'}))^{(n-1)}$.

To prove  Claim, suppose that $r_{jj'}(x)\in \partial(N_{1/2}(\hat U_{j'})$.
Then there is a
cubic $n$-dimensional symplex $K$ in $\hat U_{j'}$ such that $r_{jj'}(x)\in
\partial(N_{1/2}(K)$. The symplex $K$ is of the form $K=\prod_{i=0}^nK_i$, where
$K_i$ is a 1-dimensional symplex in $T^i_{j'}$  and there is $l\in
\{0,\dots,n\}$
such that $K_l$ is a singleton, $K_l=\{a_l\}$. Then $\partial
N_{1/2}(K)=(\cup_{k=0}^n R_k)\cap M_{j'}$, where $$R_k=
\left(\prod_{p\in\{0,\dots,n\}\setminus\{k\}}\overline{N_{1/2}(K_p)}\right)
\times \partial N_{1/2}(K_k).$$

 Then $r_{jj'}^{-1}(\partial N_{1/2}(K_k))=\partial N_{1/2}(K_k)$
and we obtain $$r_{jj'}^{-1}(R_k)\cap
M_{j'}=\left(\prod_{p\in\{0,\dots,n\}\setminus\{k\}}(r_{jj'}^p)^{-1}
(\overline{N_{1/2}(K_p)}) \right)\times \partial N_{1/2}(K_k)\cap M_{j'}.$$

For every $x=(x_0,\dots, x_n)\in r_{jj'}^{-1}(R_k)\cap
M_{j'}$ there exists $p\in\{0,\dots,n\}\setminus\{k\}$ such that
$x_p\in\partial(T^p_{j-1}\setminus T^p_j)$. Then
$d(x,(x_1,\dots,x_{k-1},y_k,x_{k+1},\dots,x_n))=1/2$ for some $y_k\in \partial
K_k$ (the boundary is considered in $T^k$; note that then $y_k\in
r^k_{j'}(\partial T^k_{j'-1}\setminus T^k_{j'})$) and
$(x_1,\dots,x_{l-1},a_l,x_{l+1},\dots,x_n)\in N_{3/4}(r_{jj'}^{-1}(\hat
U_{j'}))^{(n-1)}$. Claim is proved.

Denote by $\V$ a cover of the cubic symplex $M_j$ that splits into the union
$\V=\V^0\cup\dots\cup\V^n$ of uniformly $2^j$-bounded covers which are
$2D$-discrete.
Now for every $k$ and every $i\in\{0,\dots,n\}$ consider the family $\td
\V_k^i=\{N_{3/4}(V)\cap
r_{j,2k}^{-1}(\hat U_{2k})\mid V\in \V^i_k\}$ of  the set $r_{j,2k}^{-1}(\hat
U_{2k})$.
The cover $\V_k=\V_k^0\cup\dots\cup\V^n_k$ is uniformly $(2^j+3)$-bounded and
splits into the union of $D$-discrete families $\td\V_k^i$, $i\in\{0,\dots,n\}$.
Then the families $\W^i_k=\{r^{-1}_j(V)\in V\in\V^i_k\}$ are uniformly
$(2^{j+1}+3)$-bounded and $D$-discrete and, because
$S_{2k}=r_j^{-1}r_{j,2k}^{-1}(\hat U_{2k})$, we see that $\cup_{i=0}^n\W_k^i$
covers $S_{2k}$.

Applying Theorem \ref{t:BD} we conclude that $\asd S\le n-1$. Finally, note that
$S$ contains a separator between $A$ and $B$, namely, the set $\partial
N_{1/2}(\cup_{j=1}^\infty\hat U_j)$.
\end{proof}

\begin{prop}\label{p:sep} Let $C$ be a separator between asymptotically
disjoint subsets $A,B$ in a proper uniformly arcwise connected space $X$ and
the sets $C$ and $A\cup B$ are asymptotically disjoint. Then $C$ is an
asymptotic separator between $A$ and $B$ in $X$.
\end{prop}

\begin{proof} Denote by $U$ the set of points that can be connected to $A$ by
an arc in the set $X\setminus C$. Show that the sets $(\bar U)\cap X$ and $B$
are asymptotically disjoint. Indeed, assuming the opposite we obtain
sequences $(a_i)$ and $(b_i)$ in $U$ and $B$ respectively such
$d(a_i,x_0)\to\infty$, $d(b_i,x_0)\to\infty$ (here $x_0$ is an arbitrary base poin in $X$) and
that there
exists $R>0$ for which $d(a_i,b_i)<R$, $i\in\N$. Since $X$ is uniformly
arcwise connected, there exists $S>0$ such that, for every $i$, the points
$a_i$ and $b_i$ can be connected by an arc of diameter $\le S$.  There exists
natural $i_0$ such that for every $i\ge i_0$ we have $d(b_i,C)\ge R+S$. Then obviously
there exists an arc in $X\setminus C$ that connects $a_i$ and $b_i$, $i\ge i_0$.
 Therefore, $b_i\in U$ and we obtain a contradiction.

Thus, $B'\subset (X\setminus U)'$. In order to prove that $C$ is an asymptotic
separator between $A$ and $B$ it remains to prove that
$\bar U\cap\bar B\cap\nu X\subset C'$. Indeed, suppose the opposite, i.e. that
there exists a point $y\in \bar U\cap\bar B\cap\nu X\setminus C'$. Then there
exist disjoint neighborhoods $V$ of $y$  and  $W$ of
$\bar C$ in $\bar X$ and obviously $\overline{U\cap V} \cap \bar B\neq\emptyset$.
Therefore, there exist sequences $(a_i)$ in $U\cap V$ and $(b_i)$ in $B$ such that
$\lim_{i\to\infty}d(a_i,x_0)=\infty$, $\lim_{i\to\infty}d(b_i,x_0)=\infty$, and
the sequence $(d(a_i,b_i))$ is bounded. Arguing similarly as above we conclude that then there
exists a sequence $(c_i)$ in $C$ for which the sequence $(d(a_i,c_i))$ is bounded.
This contradicts to our choice of $V$.

\end{proof}

\begin{lemma}\label{l:uac} The space $ M(T^0,\dots,T^n)$ is  uniformly arcwise connected.
\end{lemma}
\begin{proof}
Note that for every $l$ and every $z\in r_l(M(T^0,\dots,T^n))\setminus
r_{l+1}(M(T^0,\dots,T^n))$ there exists a path of diameter $\le2^{l+1}$ in  $
M(T^0,\dots,T^n)$ connecting $z$ with $r_l(z)$.

Let $x=(x_0,\dots,x_n),y=(y_0,\dots,y_n)\in  M(T^0,\dots,T^n)$.
Let $j$ be the minimal natural number such that $d(x,y)<2^{j}$. Then
$r_{j}(x)=r_{j}(y)$. Let $J$ be a path connecting $x$ and $y$ that consists of
paths of diameter $\le2^{l+1}$ connecting $r_{l-1}(x)$ with $r_l(x)$ and of
paths of diameter $\le2^{l+1}$ connecting $r_{l-1}(y)$ with $r_l(y)$, $1\le l\le
j$ (recall that $r_0$ is the identity map). Then \begin{align*}
\diam(J)\le&\sum_{l=1}^jd(r_{l-1}(x),r_l(x))+\sum_{l=1}^jd(r_{l-1}(y),r_l(y)) \\
\le & 2\sum_{l=1}^j2^{l+1}\le2^{j+3}\le 16d(x,y).\end{align*}

\end{proof}

Now we are going to prove Theorem \ref{t:ind}. It is proved in \cite{D2} that
$\Ind
X\ge {\mathrm{Ind}} \nu X$. This result together with the classical theorem on
the comparison of ${\mathrm{Ind}}$ and ${\mathrm{dim}}$ implies the inequality
$\Ind X\ge {\mathrm{dim}}\nu X$. Suppose $\asd X<\infty$. Then $\asd
X= {\mathrm{dim}}\nu X$ (see \cite{DKU}), and we obtain $\Ind X\ge\asd X$ for
every $X$ with $\asd X<\infty$.

To prove the opposite inequality we apply induction on $\asd X$. It is known
(see \cite{D2}) that the properties $\asd X=0$ and $\Ind X=0$ are equivalent.
Assume that the inequality  $\Ind X\le\asd X$ is proved for every $X$ with $\asd
X\le n-1$.

Suppose that $\asd X\le n$. Then by Theorem \ref{t:univ} there is a large scale
embedding of $X$ into the space $M(T^0,\dots,T^n)$, for some regular locally
finite trees $T^0,\dots,T^n$. Since the image of any space under a large scale
embedding is coarse equivalent to this space and the dimension $\Ind$ is a
coarse invariant (see \cite[Proposition 1]{D2}),  we may assume that $X$ is a
subspace of  $M(T^0,\dots,T^n)$. Let $A,B$ be asymptotically disjoint subsets of
$X$. Without loss of generality we assume $A,B$ to be closed and disjoint. By
Proposition \ref{p:assep}, there exists a separator $C$ between $A$ and $B$ with $\asd
C\le n-1$. Since, by Lemma \ref{l:uac},  the space  $M(T^0,\dots,T^n)$ is uniformly
arcwise connected, it follows from Proposition \ref{p:sep} that $C$ is also an
asymptotic
separator   between $A$ and $B$ in   $M(T^0,\dots,T^n)$.

By  Lemma \ref{l:ass}, there exists  an asymptotic separator $\td C$  between
$A$ and $B$ in $X$ with $\asd \td C\le n-1$. Applying the induction assumption
we see that $\Ind \td C\le n-1$.

Since for every asymptotically disjoint subsets in $X$ there exists an
asymptotic separator between them whose  dimension $\Ind$ does not exceed $n-1$,
we conclude that $\Ind X\le n$.

The equality $\Ind X=\asd X$ is therefore proven for all $X$ with $\asd
X<\infty$.

One can similarly prove the following result.
\begin{thm} For all proper metric spaces $X$ with $\asd
X<\infty$ we have $\ind X=\asd X$.
\end{thm}

We finally obtain that the dimensions $\asd$, $\Ind$, and $\ind$ coincide in the
class of proper metric spaces that are finite dimensional with respect to
$\asd$. This fact can be regarded as a counterpart in the asymptotic topology of
the classical result on coincidence of the covering dimension, the small
inductive dimension, and the large inductive
dimension (see \cite{E}).

\section{Open problems}

The coincidence of the dimension functions $\asd$ and $\Ind$ is proved under the assumption
of finiteness of $\asd$. This leads to the following natural question.
\begin{que} Is there a proper metric space $X$ with $\asd X=\infty$ and $\Ind X<\infty$?
\end{que}

A similar question can be formulated for $\mathrm{ind}$.

 A closed subset $C$ of a topological space $X$ is a {\it cut\/}
between disjoint subsets $A, B\subset X$ if every continuum (compact connected
space) $T\subset X$ that intersect both $A$ and $B$ also intersects $C$.

The notion of the large  inductive dimension can be based on the notion of cut
instead of separator. It turnes out (see \cite{FLS}) that the obtained dimension
(it was defined by Brouwer who called it Dimensiongrad)  coincides with the classical
large inductive dimension in the class of separable metrizable spaces).

\begin{que} Does the large asymptotic inductive dimension defined on the
base of asymptotic cut coincides with $\Ind$?
\end{que}

\end{document}